\documentclass[12pt]{article}

\usepackage{graphicx,tikz,color,url}
\usepackage{amsmath,amssymb,latexsym}
\usepackage{pgfplots}
\usepackage{cite}
\usepackage{mathrsfs}
\usetikzlibrary{arrows}

\newtheorem{theorem}{Theorem}[section]

\newtheorem{lemma}[theorem]{Lemma}
\newtheorem{proposition}[theorem]{Proposition}
\newtheorem{conjecture}[theorem]{Conjecture}

\newcommand{\proof}{\noindent{\bf Proof.\ }}
\newcommand{\qed}{\hfill $\square$ \bigskip}

\newcommand{\gp}{{\rm gp}}

\newcommand{\sr}{{}_{\rm SR}}

\newcommand{\gpt}{\gp_{\rm t}}
\newcommand{\gpd}{\gp_{\rm d}}
\newcommand{\gpo}{\gp_{\rm o}}

\begin{document}

\title{On the variety of general position problems under vertex and edge removal}

\author{Jing Tian$^{a,}$\thanks{\texttt{jingtian526@126.com}}
\and 
Pakanun Dokyeesun$^{b,}$\thanks{\texttt{papakanun@gmail.com}}
\and Sandi Klav\v zar $^{c,b,d,}$\thanks{\texttt{sandi.klavzar@fmf.uni-lj.si}}
}
\maketitle

\begin{center}
$^a$ School of Science, Zhejiang University of Science and Technology, Hangzhou, Zhejiang 310023, PR China\\

$^b$ Institute of Mathematics, Physics and Mechanics, Ljubljana, Slovenia\\

\medskip

$^c$ Faculty of Mathematics and Physics, University of Ljubljana, Slovenia\\
\medskip

$^d$ Faculty of Natural Sciences and Mathematics, University of Maribor, Slovenia\\
\medskip

\medskip

\end{center}

\begin{abstract}
Let ${\rm gp}_{\rm t}(G)$, ${\rm gp}_{\rm o}(G)$, and ${\rm gp}_{\rm d}(G)$ be the total, the outer, and the dual general position number of a graph $G$, respectively. This paper investigates how removing a vertex or removing an edge affects these graph invariants. It is proved that if $x$ is not a cut vertex, then ${\rm gp}_{\rm t}(G) -1 \le {\rm gp}_{\rm t}(G-x) \le {\rm gp}_{\rm t}(G) + {\rm deg}_G(x)$. On the other hand, ${\rm gp}_{\rm o}(G-x)$ and ${\rm gp}_{\rm d}(G-x)$ can be respectively arbitrarily larger/smaller than ${\rm gp}_{\rm o}(G)$ and ${\rm gp}_{\rm d}(G)$. On the positive side, it is proved that if $x$ lies in some ${\rm gp}_{\rm o}$-set, then ${\rm gp}_{\rm o}(G)-1 \le {\rm gp}_{\rm o}(G-x)$, and that if $x$ is not a cut vertex and lies in some ${\rm gp}_{\rm d}$-set of $G$, then $ {\rm gp}_{\rm d}(G)-1 \le {\rm gp}_{\rm d}(G-x)$. For the edge removal, it is proved that (i) ${\rm gp}_{\rm t}(G) -|S(G)_{e}| \le {\rm gp}_{\rm t}(G-e) \le {\rm gp}_{\rm t}(G) +2$, where $S(G)_{e}$ is the set of simplicial vertices adjacent to both endvertices of $e$, (ii) ${\rm gp}_{\rm o}(G)/2\le {\rm gp}_{\rm o}(G-e)\leq\ 2{\rm gp}_{\rm o}(G)$, and (iii) that ${\rm gp}_{\rm d}(G) - {\rm gp}_{\rm d}(G-e)$ can be arbitrarily large. All bounds are demonstrated to be sharp. 
\end{abstract}


\medskip\noindent
{\bf AMS Subj.\ Class.\ (2020)}:  05C12, 05C69

\section{Introduction}

The concept of a general position set was introduced into graph theory independently in two papers, namely in~\cite{chandran-2016, manuel-2018a}, but in the case of hypercubes, implicitly also much earlier in~\cite{Korner-1995}. It was the article~\cite{manuel-2018a} that has stimulated a great deal of interest in this concept. The survey~\cite{survey} provides a systematic review of the results, in addition we point to the following selected recent publications~\cite{Araujo-2025, DiStefano-2025, irsic-2024, klavzar-2026, KorzeVesel-2025, Kruft-2025, Roy-2025, ThaChaTuiThoSteErs-2024, Thomas-2024a, welton-2025}. 

For a general position set, that is, a set in which no three vertices lie on a common shortest path, one requires that the vertices of the set are in general position. This requirement can be generalized by requiring that the vertices of some set lie in general position with respect to the selected set. Using this idea, and considering the natural selections for sets to be considered, in~\cite{tian-2025} a variety of general position sets was defined as follows. First, for a given graph $G = (V (G), E(G))$ and a set of vertices $Z \subseteq V (G)$, two vertices $x, y \in V(G)$ are said to be {\em $Z$-positionable} if any shortest $x,y$-path intersects $Z$ at most in $x$ and $y$, that is, no inner vertex of the path lies in $Z$. Second, if $\overline{Z} = V(G)\setminus Z$, then the set $Z$ is
\begin{enumerate}
\item[(i)] a {\em general position set}, if every $u, v \in Z$ are $Z$-positionable;
\item[(ii)] a {\em total general position set}, if every $u, v \in V(G)$ are $Z$-positionable;
\item[(iii)] an {\em outer general position set}, if every $u, v \in Z$ are $Z$-positionable, and every $u \in Z$, $v \in \overline{Z}$ are $Z$-positionable; and
\item[(iv)]  a {\em dual general position set}, if every $u, v \in Z$ are $Z$-positionable, and every $u, v \in \overline{Z}$ are $Z$-positionable.
\end{enumerate} 
Typically, we are interested in the cardinalities of largest such sets, which are respectively denoted by $\gp(G)$, $\gpt(G)$, $\gpo(G)$, and $\gpd(G)$, and called the {\em general position number}, the {\em total general position number}, the {\em outer general position number}, and the {\em dual general position number} of $G$, respectively. If $\tau \in \{\gp, \gpt, \gpo, \gpd\}$, then a set $Z$ is a {\em $\tau$-set} if it has the corresponding property and $|Z| = \tau(G)$. Note that, by definition, $\gpt(G)\le \gpo(G)\le \gp(G)$ and $\gpt(G)\le \gpd(G)\le \gp(G)$. After being  introduced in~\cite{tian-2025}, this variety has been also investigated on strong and lexicographic products of graphs~\cite{dokyeesun-2026}. 

The paper~\cite{dokyeesun-2025} investigates how removing a vertex or an edge affects the general position number of a graph. It is proved that ${\rm gp}(G-x)\leq 2{\rm gp}(G)$ holds for any vertex $x$ of a connected graph $G$ and that if $x$ lies in some ${\rm gp}$-set of $G$, then ${\rm gp}(G) - 1 \le {\rm gp}(G-x)$. On the other hand, ${\rm gp}(G-x)$ can be much larger than ${\rm gp}(G)$ also when $G-x$ is connected. For the edge removal it is proved that $\gp(G)/2\le \gp(G-e)\leq\ 2\gp(G)$ holds for any edge $e$ of $G$. In this paper we continue this line of research by investigating how removing a vertex or an edge affects the other three general position numbers. 

In the next section we give definitions and recall results needed in the rest of the paper. In Section~\ref{sec:vertex} we consider the vertex removal. For the total general position sets we prove that if $x$ is not a cut vertex, then $\gpt(G) -1 \le \gpt(G-x) \le \gpt(G) + \deg_G(x)$. Section~\ref{sec:outer-vertex} addresses the outer general position number of a graph when a vertex is removed and shows that it cannot be bounded in terms of the outer general position number of the original graph. 
On the other hand, if $x$ lies in some $\gpo$-set, then $\gpo(G)-1 \le \gpo(G-x)$. Moreover, if $x$ is a simplicial vertex, then $\gpo(G-x) \le \gpo(G)+\deg_G(x)-1$. Section~\ref{sec:dual-vertex} then deals with the dual general position number. Similarly, as for the outer general position number, $\gpd(G-x)$ can be arbitrarily larger/smaller than $\gpd(G)$. On the other hand, if $x$ is not a cut vertex and lies in some $\gpd$-set of $G$, then $ \gpd(G)-1 \le \gpd(G-x)$. In the final section we consider the edge removal operation. For the total and the outer general position number we prove sharp lower and upper bounds, while for the dual general position number we show that the difference $\gpd(G) - \gpd(G-e)$ can be arbitrarily large. All bounds proved in this paper are demonstrated to be sharp.

\section{Preliminaries}

We consider simple and connected graphs $G=(V(G), E(G))$. For a positive integer $k$, the set $\{1,\ldots, k\}$ is denoted by $[k]$.

Let $u\in V(G)$. Then $N_G(u)$ denotes the set of neighbors of $u$ in $G$ and $N_G[u] = N_G(u)\cup \{u\}$. The {\em degree} $\deg_G(u)$ of $u$ is $\deg_G(u) = |N_G(u)|$. If $X\subseteq V(G)$, then the subgraph of $G$ induced by $X$ is denoted by $G[X]$. The graph $G-u$ is obtained by deleting $u$ and all incident edges from $G$, that is, $G-u$ is the subgraph $G[V(G)\setminus \{u\}]$.
For an edge $e \in E(G)$, the graph $G-e$ is obtained by deleting the edge $e$ from $G$.
A vertex $u$ of $G$ is {\em simplicial} if $N_G(u)$ induces a complete subgraph. The set of all simplicial vertices of $G$ will be denoted by $S(G)$ and the cardinality of $S(G)$ by $s(G)$. Moreover, $\omega(G)$ and $\alpha(G)$ stand for the clique number and the independence number of $G$. 

The {\em distance} $d_G(u,v)$ between vertices $u$ and $v$ of $G$ is the number of edges on a shortest $u,v$-path.  
A set $X\subseteq V(G)$ is {\em convex} if for any vertices $u,v\in X$, any shortest $u,v$-path contains only vertices from $X$. By abuse of language, we will also say that a subgraph $H$ of $G$ is convex if $V(H)$ is convex. The vertex $u$ of $G$ is \emph{maximally distant} from a vertex $v\in V(G)$ if every $w\in N_G(u)$ satisfies $d_G(v,w)\le d_G(u,v)$. If $u$ is maximally distant from $v$, and $v$ is also maximally distant from $u$, then $u$ and $v$ are \emph{mutually maximally distant}, MMD for short. Note that true twins, that is, vertices $u$ and $v$ with $N_G[u] = N_G[v]$, are MMD. The {\em strong resolving graph} $G\sr$ of $G$ has $V(G)$ as the vertex set, two vertices being adjacent in $G\sr$ if they are MMD in $G$. This concept was introduced in~\cite{oellermann-2007} for the purpose of better understanding the strong metric dimension of graphs.

The first known result which we need later on describes general position sets in an arbitrary graph. To state it, some more definitions are required. If ${\cal P} = \{X_1, \ldots, X_t\}$ is a partition of $X\subseteq V(G)$, then ${\cal P}$ is \emph{distance-constant} if for any $i,j\in [t]$, $i\ne j$,
there exists a constant $p_{ij}$, such that $d_G(x,y) = p_{ij}$ for every $x\in X_i$, $y\in X_j$. If so, we set $d_G(X_i,X_j) = p_{ij}$. A distance-constant partition ${\cal P}$ is {\em in-transitive} if $p_{ik} \ne p_{ij} + p_{jk}$ holds for $i,j,k\in [t]$.

\begin{theorem} {\rm \cite[Theorem 3.1]{anand-2019}}
\label{thm:gpsets}
Let $G$ be a graph. Then $X\subseteq V(G)$ is a general position set if and only if the components of $G[X]$ are complete subgraphs, the vertices of which form an in-transitive, distance-constant partition of $X$.
\end{theorem}	

The following theorem is the most important for us since it contains characterizations of the three variants of general position sets. They were respectively proved in~\cite[Theorems~2.1, 2.3, 3.1]{tian-2025}. 

\begin{theorem} \label{thm:all}
If $G$ is a graph and $X\subseteq V(G)$, then the following hold.
\begin{enumerate}
\item[(i)] $X$ is a total general position set of $G$ if and only if $X\subseteq S(G)$. Consequently, $\gpt(G)=s(G)$.
\item[(ii)] If $|X|\geq 2$, then $X$ is an outer general position set of $G$ if and only if each pair of vertices from $X$ is MMD. Consequently, $\gpo(G) = \omega(G\sr)$.
\item[(iii)] If $X$ is a general position set of $G$, then $X$ is a dual general position set if and only if $G-X$ is convex in $G$.
\end{enumerate}
\end{theorem} 

We conclude the preliminaries with the following straightforward, but very useful observation. 

\begin{lemma}
\label{lem:s(G)-is-lower-bound}
If $G$ is a graph and $\tau\in \{\gp, \gpo, \gpd, \gpt\}$, then $\tau(G)\ge s(G)$. Moreover, the equality holds if $G$ is a block graph. 
\end{lemma}

\section{Vertex removal}
\label{sec:vertex}

Before we turn our attention to  vertex-deleted subgraphs in specific variants of general position sets, we consider the following example which is useful to understand each of the three invariants investigated. Let $G_n$, $n\ge 2$, be the graph, and $x$ its vertex as shown in Fig.~\ref{fig:gpo-vertex}. 

\begin{figure}[ht!]
        \centering
\begin{tikzpicture}[line cap=round,line join=round,>=triangle 45,x=1.5cm,y=1.5cm]
\draw [line width=1pt](2.,4.5)-- (2.,3.);
\draw [line width=1pt] (2.,3.)-- (4.,3.);
\draw [line width=1pt] (4.,3.)-- (6.,3.);
\draw [line width=1pt] (6.,3.)-- (6.,4.5);
\draw [line width=1pt] (2.,4.5)-- (3.,5.);
\draw [line width=1pt] (2.,4.5)-- (3.,4.);
\draw [line width=1pt] (2.,4.5)-- (3.,3.5);
\draw [line width=1pt] (4.,4.5)-- (3.,4.);
\draw [line width=1pt] (4.,4.5)-- (3.,3.5);
\draw [line width=1pt] (3.,5.)-- (4.,4.5);
\draw [line width=1pt] (4.,4.5)-- (5.,5.5);
\draw [line width=1pt] (5.,5.5)-- (6.,4.5);
\draw [line width=1pt] (2.,4.5)-- (1.,3.5);
\draw [line width=1pt] (6.,4.5)-- (7.,3.5);
\draw [line width=1pt] (2.,4.5)-- (1.,4.);
\draw [line width=1pt] (6.,4.5)-- (7.,4.);
\draw [line width=1pt] (2.,4.5)-- (1.,5.);
\draw [line width=1pt] (6.,4.5)-- (7.,5.);
\draw [line width=1pt] (2.,4.5)-- (1.,5.5);
\draw [line width=1pt] (6.,4.5)-- (7.,5.5);
\draw [line width=1pt] (2.,4.5)-- (3.,5.5);
\draw [line width=1pt] (3.,5.5)-- (4.,4.5);
\draw [line width=1pt] (4.,4.5)-- (5.,5);
\draw [line width=1pt] (5.,5)-- (6.,4.5);
\draw [line width=1pt] (4.,4.5)-- (5.,4.);
\draw [line width=1pt] (4.,4.5)-- (5.,3.5);
\draw [line width=1pt] (6.,4.5)-- (5.,4.);
\draw [line width=1pt] (6.,4.5)-- (5.,3.5);
\draw (3.8,4.2) node[anchor=north west] {$x$};
\draw (1.7,2.8) node[anchor=north west] {$u$};
\draw (3.8,2.8) node[anchor=north west] {$v$};
\draw (5.7,2.8) node[anchor=north west] {$w$};
\draw (4.85,4.85) node[anchor=north west] {$\vdots$};
\draw (2.85,4.85) node[anchor=north west] {$\vdots$};
\draw (6.85,4.85) node[anchor=north west] {$\vdots$};
\draw (0.85,4.85) node[anchor=north west] {$\vdots$};
\draw[decoration={brace,mirror,raise=5pt},decorate]
  (7.3,3.5) -- node[right=6pt] {$n$} (7.3,5.5);
  \draw[decoration={brace,raise=5pt},decorate]
  (0.7,3.5) -- node[left=6pt] {$n$} (0.7,5.5);
\begin{scriptsize}
\draw [fill=white] (2.,4.5) circle (2pt);
\draw [fill=white] (6.,4.5) circle (2pt);
\draw [fill=white] (2.,3.) circle (2pt);
\draw [fill=white] (4.,3.) circle (2pt);
\draw [fill=white] (6.,3.) circle (2pt);
\draw [fill=white] (3.,5.) circle (2pt);
\draw [fill=white] (3.,4) circle (2pt);
\draw [fill=white] (3.,3.5) circle (2pt);
\draw [fill=white] (5.,4) circle (2pt);
\draw [fill=white] (5.,3.5) circle (2pt);
\draw [fill=white] (4.,4.5) circle (2pt);
\draw [fill=white] (5.,5) circle (2pt);
\draw [fill=white] (1.,3.5) circle (2pt);
\draw [fill=white] (7.,3.5) circle (2pt);
\draw [fill=white] (1.,4.) circle (2pt);
\draw [fill=white] (7.,4.) circle (2pt);
\draw [fill=white] (1.,5.) circle (2pt);
\draw [fill=white] (7.,5.) circle (2pt);
\draw [fill=white] (1.,5.5) circle (2pt);
\draw [fill=white] (7.,5.5) circle (2pt);
\draw [fill=white] (3.,5.5) circle (2pt);
\draw [fill=white] (5.,5.5) circle (2pt);

\draw[rounded corners] (2.7, 3.3) rectangle (3.3, 6.0);
\draw (2.8,6.0) node[anchor=north west] {$X_1$};
\draw[rounded corners] (4.7, 3.3) rectangle (5.3, 6.0);
\draw (4.8,6.0) node[anchor=north west] {$X_2$};
\draw[rounded corners] (6.7, 3.3) rectangle (7.3, 6.0);
\draw (6.8,6.0) node[anchor=north west] {$Y$};

\end{scriptsize}
\end{tikzpicture}
        \caption{Graph $G_n$ with $|X_1|=|X_2|=n$.}
        \label{fig:gpo-vertex}
    \end{figure}

\begin{proposition}
\label{prop:graph-Gn}
If $n\ge 2$ and $\tau\in \{\gpo, \gpd, \gpt\}$, then $$\tau(G_n)=2n\quad {\rm and}\quad  \tau(G_n-x)=4n\,.$$    
\end{proposition}

\proof
Since $s(G_n) = 2n$, Theorem~\ref{thm:all}(i) yields $\gpt(G_n) = 2n$. Observe that two vertices of $G_n$ are MMD if and only if they are leaves, or they both belong to $X_1$, or they both belong to $X_2$. Hence the strong resolving graph $(G_n)\sr$ has three non-trivial components, one isomorphic to $K_{2n}$, and the other two each isomorphic to $K_n$. Hence $\gpo(G_n) = 2n$ by Theorem~\ref{thm:all}(ii). By Lemma~\ref{lem:s(G)-is-lower-bound} we also have $\gpd(G_n)\geq 2n$. To prove the reverse inequality, we will apply Theorem~\ref{thm:all}(iii). For this sake we first get by a case analysis that $\gp(G_n) = 2n+1$. By symmetry, the $\gp$-sets to be considered are $X_1\cup X_2 \cup\{v\}$, $X_2\cup Y \cup\{u\}$, $X_2\cup Y \cup\{v\}$, $X_2\cup Y \cup\{w\}$, $X_1\cup Y \cup\{u\}$, $X_1\cup Y \cup\{v\}$, and $X_1\cup Y \cup\{w\}$. If $Z$ is any of these sets, then $G_n - Z$ is not connected, hence clearly not convex. Theorem~\ref{thm:all}(iii) thus implies that none of these sets is a dual general position set. It follows that $\gpd(G_n)\le 2n$ and we can conclude that $\gpd(G_n) = 2n$. 

Since the graph $G_n - x$ is a tree, Lemma~\ref{lem:s(G)-is-lower-bound} gives
$$\gpt(G_n-x) = \gpo(G_n-x) = \gpd(G_n-x) = s(G_n-x) = 4n\,,$$ 
and we are done. 
\qed

\subsection{Total general position sets}
\label{sec:total}

In this subsection we show that the total general position set of $G-x$ can be bounded below and above when $x$ is not a cut vertex of a graph $G$.  

\begin{theorem}
\label{thm:gpt-minus-vertex}
If $x$ is not a cut vertex of a graph $G$, then 
$$\gpt(G) -1 \le \gpt(G-x) \le \gpt(G) + \deg_G(x)\,.$$
Moreover, if  $x$ is a simplicial vertex, then 
$$\gpt(G-x) \le \gpt(G) + \deg_G(x) - 1\,.$$ 
In addition, all three bounds are sharp. 
\end{theorem}

\proof
The bounds clearly hold if $G$ is complete, hence assume in the rest of the proof that $G$ is not complete. 

Let $X$ be a $\gpt$-set of $G$. Then $X = S(G)$ by Theorem~\ref{thm:all}(i). 

To prove the lower bound, it suffices to show that $X\setminus\{x\} \subseteq S(G-x)$. Suppose to the contrary that $y \in X\setminus\{x\}$ is not simplicial in $G-x$. Then there exist $z,w \in N_{G-x}(y)$ such that $z$ and $w$ are not adjacent in $G-x$. The vertices $z$ and $w$ are also not adjacent in $G$, but then $y$ is not simplicial in $G$, a contradiction. Hence $X\setminus\{x\} \subseteq S(G-x)$ and by applying Theorem~\ref{thm:all}(i) we get $\gpt(G) -1 \le |X\setminus\{x\}| \le \gpt(G-x)$. This proves the lower bound. 

For the upper bound observe that if $y\in V(G)\setminus N_G[x]$ is a simplicial vertex of $G-x$, then $y$ is a simplicial vertex of $G$. Using Theorem~\ref{thm:all}(i) once more we can argue as follows:
$$\gpt(G-x) = s(G-x) \le s(G) + \deg_G(x) = \gpt(G) + \deg_G(x)\,.$$
This proves the upper bound in the general case.
Assume now that $x$ is a simplicial vertex. Since $G$ is not complete, at least one vertex in $N_G(x)$ is not simplicial in $G-x$, hence  
$$\gpt(G-x) = s(G-x) \le s(G) + (\deg_G(x) - 1) = \gpt(G) + \deg_G(x) - 1\,.$$

To show that the lower bound is sharp, consider the stars $K_{1,n}, n \ge 3$. Clearly, $s(K_{1,n}) = n$ and then $s(K_{1,n}-x) = n-1$ for any leaf $x$ of $K_{1,n}$. 
To demonstrate that the upper bound is tight, consider complete bipartite graphs $K_{2,n}$, $n\ge 3$. If $x$ is a vertex of $K_{2,n}$ of degree $n$, then we have $\gpt(K_{2,n}) = 0$ and $\gpt(K_{2,n}-x) = n = \gpt(K_{2,n})+\deg_{K_{2,n}}(x)$. Finally, to show that the upper bound is sharp in the case when $x$ is a simplicial vertex, consider the edge deleted complete graph $K_n-e$ and let $x_1$ and $x_2$ be the vertices of $K_n-e$ of degree $n-2$. Then $S(K_n-e) = \{x,y\}$, so that $\gpt(K_n-e) = 2$. On the other hand, $(K_n-e)-x_1\cong K_{n-1}$. Hence $\gpt((K_n-e)-x_1) = n-1 = \gpt(K_n-e) + \deg_{K_n-e}(x_1) - 1$.
\qed

Note that in view of Theorem~\ref{thm:all}(i), we can rephrase Theorem~\ref{thm:gpt-minus-vertex} by saying that if $x$ is not a cut vertex of a graph $G$, then 
$$s(G) -1 \le s(G-x) \le s(G) + \deg_G(x)\,,$$
and if $x$ is a simplicial vertex, then 
$$s(G-x) \le s(G) + \deg_G(x) - 1\,.$$

\subsection{Outer general position sets}
\label{sec:outer-vertex}

In this subsection we show that anything can happen with outer general position sets when a vertex of a graph $G$ is removed. 

To demonstrate that $\gpo(G-x)$ can be arbitrarily larger than $\gpo(G)$, consider the graph $G_n$, $n\ge 2$, from Fig.~\ref{fig:gpo-vertex}. By Proposition~\ref{prop:graph-Gn}, $\gpo(G_n) = 2n$ and $\gpo(G_n-x) = 4n$.

To show that $\gpo(G-x)$ can also be arbitrarily smaller than $\gpo(G)$, consider the {\em fan graph} $F_n$, $n\ge 3$,  of order $n+1$, which is obtained from $P_n$ and a vertex adjacent to all the vertices of $P_n$. Let $x$ be the vertex of $F_n$ of degree $n$. In~\cite[Proposition 3.2]{dokyeesun-2026} it is proved that if $G$ is a diameter $2$ graph with no true twins, then $\gpo(G) = \alpha(G)$. Since $F_n-x\cong P_n$, we thus have 
$$\gpo(F_n) = \lceil n/2\rceil\quad {\rm and}\quad \gpo(F_n-x) = 2\,.$$

On the other hand, in the case when a vertex $x$ lies in some $\gpo$-set, we can bound $\gpo(G-x)$ from below as follows. 

\begin{theorem}
\label{thm:outer-vertex}
    If $x$ is not a cut vertex and lies in some $\gpo$-set of a graph $G$, then $ \gpo(G)-1 \le \gpo(G-x)$ and the bound is sharp.
\end{theorem}

\proof
Let $X$ be a $\gpo$-set of $G$ such as $x \in X$. Then $|X| = \gpo(G)$ and by Theorem~\ref{thm:all}(ii), every pair of vertices in $X$ are MMD in $G$.

To show the bound, it suffices to prove the assertion that $X\setminus\{x\}$ is an outer general position set of $G-x$. Given Theorem~\ref{thm:all}(ii) this is equivalent to proving that every two vertices in $X\setminus \{x\}$ are MMD in $G-x$. Suppose to the contrary that there exist $u, v \in X\setminus \{x\}$ such that $u$ and $v$ are not MMD in $G-x$. By symmetry we may assume that there exists $u' \in N_{G-x}(u)$ such that $d_{G-x}(u',v) = d_{G-x}(u,v) + 1$. Since $u$ and $v$ are MMD in $G$, we have $d_{G}(u',v) \le d_{G}(u,v)$. This implies that $x$ lies on some shortest $u',v$-path in $G$. But then $v$ and $x$ are not MMD in $G$, a contradiction. This contradiction implies that every pair of vertices in $X\setminus \{x\}$ are MMD in $G-x$. We can conclude that  $\gpo(G-x) \ge |X\setminus \{x\}| = \gpo(G) -1$.

For sharpness of the bound, consider the stars $K_{1,n}$, $n \ge 3$. Then $\gpo(K_{1,n}) = n$ and $\gpo(K_{1,n}-x) = n-1$ for any leaf $x$ of $K_{1,n}$.
\qed

By Theorem~\ref{thm:outer-vertex} we can thus bound $\gpo(G-x)$ from below when $x$ is not a cut vertex and lies in some $\gpo$-set of $G$. As a possible general upper bound on $ \gpo(G-x)$ we pose: 

\begin{conjecture}
\label{con:gpo}
If $x$ is not a cut vertex of a graph $G$, then $$\gpo(G-x) \le \gpo(G) + \deg_G(x)\,.$$ 
\end{conjecture}

If Conjecture~\ref{con:gpo} holds true, then it is sharp as justified by Proposition~\ref{prop:graph-Gn}. For the case when $x$ is a simplicial vertex, the conjecture holds true, and even more is true: 

\begin{proposition}
\label{prop:outer-simplicial-upper}
If $x$ is a simplicial vertex of a graph $G$, then 
$$ \gpo(G-x) \le \gpo(G)+\deg_G(x)-1$$ 
and the bound is sharp.
\end{proposition}

\proof
Let $u$ and $v$ be arbitrary vertices of $G-x$. Since $x$ is a simplicial vertex of $G$, no shortest $u,v$-path in $G$ contains $x$ which in turn implies that a $u,v$-path is shortest in $G$ if and only if it is shortest in $G-x$. Consequently, if $u,v\in V(G)\setminus N_G[x]$, then $u, v$ are MMD in $G$ if and only if they are MMD in $G-x$. 

Let $\gpo(G) = k$ and let $X$ be a largest set of pairwise MMD vertices of $G$. Note that $x\in X$. Let further $Y$ be a largest set of pairwise MMD vertices of $G-x$. Then by the above, $|Y\cap (V(G)\setminus N_G[x])| \le k-1$. This in turn implies that 
$$|Y| \le (k-1) + |N_G(x)| = \gpo(G) - 1 + \deg_G(x)\,,$$
which proves the bound. 

Let $G_{n,k}$, $k < \frac{n+1}{2}$, be the graph obtained from $K_n$ and another vertex $x$ which is adjacent to $k$ vertices of $K_n$. Then $\gpo(G_{n,k}) = \max \{n-k+1, k\} = n-k+1$, where the last equality holds because we have assumed that $k < \frac{n+1}{2}$. On the other hand, $G_{n,k}-x \cong K_n$ and hence $\gpo(G_{n,k}-x) = n$. 
\qed

\subsection{Dual general position sets}
\label{sec:dual-vertex}

In this subsection we first show that $\gpd(G-x)$ can be arbitrarily larger/smaller than $\gpd(G)$. After that we prove that if $x$ is not a cut vertex and lies in some $\gpd$-set of $G$, then $\gpd(G)-1 \le \gpd(G-x)$. We conclude the subsection by demonstrating that for such a vertex $x$, the value  $\gpd(G-x)$ cannot be bounded from above by $\gpd(G)$. 

Consider the fan graphs $F_n$ ($n\geq 4$) defined in Subsection 3.2, where $x$ denotes the vertex of degree $n$. It was proved in~\cite[Corollary 2.9]{tian-2023} that if $n\ge 4$, then $\gp(F_n) = \left\lfloor \frac{2(n+1)}{3}\right\rfloor$. Moreover, a gp-set of $F_n$ consists of the end-vertices of a largest set of independent edges of the path $P_n$ in $F_n$. By Theorem~\ref{thm:all}(iii) it is then straightforward to check that every largest general position set of $F_n$ is a dual general position set. Hence $\gp(F_n) = \gpd(F_n)$. Moreover, since $\gp(F_n) = \left\lfloor \frac{2(n+1)}{3}\right\rfloor$ for $n\ge 4$, and because $F_n - x\cong P_n$, we have 
$$\gpd(F_n) = \left\lfloor \frac{2(n+1)}{3}\right\rfloor \quad {\rm and}\quad \gpd(F_n-x) = 2\,.$$ 

On the other hand, consider the family of graphs defined as follows. 
The {\em wheel} $W_{n}$, $n\ge 4$, is the graph of order $n+1$ obtained from $C_n$ by adding one more vertex and connecting it to all vertices of the cycle. The {\em mushroom} $M_k$, $k\ge 4$, is the graph obtained from the disjoint union of $W_{k+4}$ and $K_k$ by adding a matching between the vertices of $K_k$ and $k$ consecutive vertices of the $(k+4)$-cycle of $W_{k+4}$. See Fig.~\ref{fig:mushroom} for $M_4$. 

\begin{figure}[ht!]
	\centering
	\begin{tikzpicture}[line cap=round,line join=round,>=triangle 45,x=1cm,y=1cm]
\draw[line width=1pt] (0,0) -- (1,0) -- (2,0) -- (3,0) -- (4,0)--(5,0)--(6,0)--(7,0);
\draw[line width=1pt] (3.5,2.5)--(0,0);
\draw[line width=1pt] (3.5,2.5)--(1,0);
\draw[line width=1pt] (3.5,2.5)--(2,0);
\draw[line width=1pt] (3.5,2.5)--(3,0);
\draw[line width=1pt] (3.5,2.5)--(4,0);
\draw[line width=1pt] (3.5,2.5)--(5,0);
\draw[line width=1pt] (3.5,2.5)--(6,0);
\draw[line width=1pt] (3.5,2.5)--(7,0);	
\draw[line width=1pt] (0,0) .. controls (0,-1) and (7,-1) .. (7,0);
\draw[line width=1pt] (2,0)--(2,-2);  
\draw[line width=1pt](3,0) --(3,-2); 
\draw[line width=1pt](4,0)--(4,-2); 
\draw[line width=1pt](5,0)--(5,-2);
\draw[line width=1pt](2,-2)--(3,-2)--(4,-2)--(5,-2);
\draw[line width=1pt] (2,-2) .. controls (2,-2.5) and (4,-2.5) .. (4,-2);
\draw[line width=1pt] (3,-2) .. controls (3,-2.5) and (5,-2.5) .. (5,-2);
\draw[line width=1pt] (2,-2) .. controls (2,-3) and (5,-3) .. (5,-2);
\draw (3.5,2.8) node {$x$};
\draw (-0.1,.25) node {$w_7$};
\draw(0.8,0.25) node {$w_8$};
\draw(1.8,0.25) node {$w_1$};
\draw(2.8,0.25) node {$w_2$};
\draw(4.2,0.25) node {$w_3$};
\draw(5.2,0.25) node {$w_4$};
\draw(6.2,0.25) node {$w_5$};
\draw(7.1,0.25) node {$w_6$};
\draw(1.75,-1.75) node {$v_1$};
\draw(2.75,-1.75) node {$v_2$};
\draw(4.25,-1.75) node {$v_3$};
\draw(5.25,-1.75) node {$v_4$};
\begin{scriptsize}
\draw [fill=white] 
(0,0) circle (2pt)
(1,0) circle (2pt)
(2,0) circle (2pt)
(3,0) circle (2pt)
(4,0) circle (2pt)
(5,0) circle (2pt)
(6,0) circle (2pt)
(7,0) circle (2pt)
(3.5,2.5) circle (2pt)
(2,-2) circle (2pt)
(3,-2) circle (2pt)
(4,-2) circle (2pt)
(5,-2) circle (2pt);
\end{scriptsize}
\end{tikzpicture}
\caption{The mushroom $M_4$.}
\label{fig:mushroom}
\end{figure}

\begin{proposition}
\label{prop:dual-mushroom graph}
If $k\geq 4$, then $\gpd(M_k) = k+2$ and $\gpd(M_k - x) = 0$.
\end{proposition}

\proof
Let $k\geq 4$, and set 
$$V(M_k) =  \{w_1,\ldots,w_{k+4}, x, v_1,\ldots, v_k\}\,,$$
where $V(W_{k+4}) = \{w_1,\ldots,w_{k+4}, x\}$ with $x$ being the center of $W_{k+4}$, $V(K_{k}) = \{v_1,\ldots, v_{k}\}$, and $v_iw_i\in E(M_k)$ for $i\in [k]$, see Fig.~\ref{fig:mushroom} again.

Let $Y = V(K_k)\cup\{w_{k+2},w_{k+3}\}$. Then it follows from Theorem~\ref{thm:gpsets} that $Y$ is a general position set of $M_k$. Moreover, $M_k - Y$ is convex, hence Theorem~\ref{thm:all}(iii) implies that $Y$ is a dual general position set of $M_k$. Thus $\gpd(M_k)\geq k+2$.

To prove that $\gpd(M_k)\leq k+2$, suppose to the contrary that there exists a dual general position set $X$ of $M_k$ with $|X| > k+2$. We first claim that $x\notin X$ by using the argument that if $x$ would belong to $X$, then $G-X$ would not be convex. So suppose that $x\in X$. Then at least one of $w_1$ and $w_4$ must belong to $X$, say $w_1\in X$. Then $w_4\notin X$. Since $w_1\in X$, exactly one of $w_2$ and $w_{k+4}$ lies in $X$.  If $w_{k+4}\in X$, then $w_2$ and $w_4$ are not $X$-positionable. And if $w_{2}\in X$, then $w_{k+4}$ and $w_4$ are not $X$-positionable. We analogously get a contradiction in the case when $w_{4}\in X$. We can conclude that $x\notin X$. We next claim that $w_i\notin X$ for each $i\in [k]$. Suppose that $w_i\in X$ for some $i\in [k]$. Since $x\notin X$, we get that $v_i\in X$ because $x,w_i,v_i$ is the unique shortest $x,v_i$-path. Further, exactly one of the neighbors of $w_{i}$ on the $(k+4)$-cycle of $W_{k+4}$ lies in $X$. But this is not possible, as $X$ is a general position set.  

From our assumption, we have $|X\cap\{w_{k+1},w_{k+2},w_{k+3},w_{k+4}\}|\geq 3$. 
If $w_{k+1}\in X$, then $w_{k+2}\in X$ because $w_k\notin X$. Otherwise, the vertices $w_k$ and $w_{k+2}$ are not $X$-positionable. 
Since $w_{k+2}$ lies on the shortest $w_{k+1},w_{k+3}$-path, it follows that $w_{k+3}\notin X$ and $w_{k+4}\in X$. If $w_{k+4}\in X$, then we must have $w_1\in X$ which is not possible because we have proved that $w_i\notin X$ for each $i\in [k]$.  This contradiction implies that $w_{k+1}\notin X$. Similarly, we get that $w_{k+4}\notin X$.  Therefore, $|X\cap \{w_{k+1},w_{k+2},w_{k+3},w_{k+4}\}|\leq 2$, which leads to $|X|\leq k+2$, the final contradiction. We can conclude that $\gpd(M_k) = k+2$. 

We now prove that $\gpd(M_k-x)=0$. Let $X$ be an arbitrary dual general position set of $M_k-x$. Then we first infer that $w_i\notin X$ for each $i\in [k]$. Indeed, if $w_i\in X$ for some $i\in [k]$, then exactly one of $w_{i-1}, w_{i+1}$ must belong to $X$, but this leads to a contradiction that $X$ is a dual general position set. Similarly we get that none of $w_{k+j}$, $j\in [4]$, belong to $X$. Finally, if some $v_i\in X$, then we get that all $v_i$, $i\in [k]$ belong to $X$, but then $w_1$ and $w_k$ are not $X$-positionable. We can conclude that $X = \emptyset$.
\qed

Proposition~\ref{prop:dual-mushroom graph} shows that $\gpd(G - x)$ can be arbitrarily smaller than $\gpd(G)$. Next, we show that $\gpd(G - x)$ can also be arbitrarily larger than $\gpd(G)$ by the following example. Let $T_k$, $k\ge 3$, be a graph obtained from the disjoint union of $K_{1,2k}$ with the central vertex $u$ and leaves $v_1,\ldots, v_{2k}$, and an isolated vertex $x$, by adding the edges $xv_i$, $i \in [k]$.

\begin{proposition}
\label{prop:dual-increase}
If $k\geq 3$, then $\gpd(T_k) = k$ and $\gpd(T_k - x) = 2k$.
\end{proposition}

\proof
Since $s(T_k) = k$,  Lemma~\ref{lem:s(G)-is-lower-bound} gives $\gpd(T_k)\geq k$.  
To prove the reverse inequality, consider an arbitrary $\gpd$-set $X$ of $T_k$. Suppose first that $x\in X$. Then $v_i, v_j\in X$, where $i,j \in [k]$, is not possible since then $v_i$ and $v_j$ are not $X$-positionable. Analogously we see that also $v_i, v_j\notin X$, where $i,j \in [k]$, cannot happen. Since $k\ge 3$, we can conclude that $x\notin X$. By an analogous reason, $u \notin X$. This in turn implies that $v_i\notin X$ for each $i\in [k]$. Hence $|X|\leq k$, and then $\gpd(T_k) = k$.

Since $s(T_k - x) = 2k$, Lemma~\ref{lem:s(G)-is-lower-bound} implies that $\gpd(T_k-x) = 2k$, and we are done.
\qed

We also have a result parallel to Theorem~\ref{thm:outer-vertex} for the dual general position number in the following. 

\begin{theorem}
\label{thm:dual-vertex}
If $x$ is not a cut vertex and lies in some $\gpd$-set of a graph $G$, then $\gpd(G)-1 \le \gpd(G-x)$ and the bound is sharp.
\end{theorem}

\proof
Let $X$ be a $\gpd$-set of $G$ such that $x \in X$. Then $|X| = \gpd(G)$  and by Theorem~\ref{thm:all}(ii), $G-X$ is convex. 
By the proof of~\cite[Proposition 3.3]{dokyeesun-2025}, the set $X\setminus\{x\}$ is a general position set of $G-x$. It suffices to show that $(G-x) -(X\setminus\{x\})$ is convex in $G-x$, because then by Theorem~\ref{thm:all}(iii) the set $X\setminus\{x\}$ is a dual general position set of $G-x$ and hence $\gpd(G-x) \ge |X\setminus \{x\}| = \gpd(G) -1$. 

Suppose to the contrary that this is not the case. There are $u,v \notin X$ such that a shortest $u,v$-path $P$ contains a vertex $w \in X\setminus \{x\}$. Since $G-X$ is convex, $d_G(u,v) < d_G(u,w)+d_G(w,v)$. Based on this, we can conclude the following: 
\begin{align*}
d_{G-x}(u,v) & = d_{G-x}(u,w)+d_{G-x}(w,v)\\
& \ge  d_G(u,w)+d_G(w,v)\\
& > d_G(u,v)\,.
\end{align*}
This in turn implies that $x$ lies on the shortest $u,v$-path in $G$ which contradicts with $G-X$ being convex. We can  conclude that  $(G-x) -(X\setminus\{x\})$ is convex in $G-x$. 

To prove the sharpness of the bound, consider a graph $K_{1,n}$, $n\geq 3$, and let $x$ be an arbitrary leaf of it. Then $\gpd(K_{1,n})=n$ and $\gpd(K_{1,n}-x) =n-1$. 
\qed

In Theorem~\ref{thm:dual-vertex} we have bounded $\gpd(G-x)$ from below by $\gpd(G)-1$ when $x$ is not a cut vertex and lies in some $\gpd$-set of $G$. We next demonstrate that for such a vertex $x$ the value  $\gpd(G-x)$ cannot be bounded from above by $\gpd(G)$. 

Consider the graph $Y_k$, which is obtained from the complete graph $K_{k+3}$ by adding two vertices $x$ and $y$, and connecting each of them them to the same $k$ vertices of $K_{k+3}$ as illustrated in Fig.~\ref{fig:Graph $Y_k$}. 

\begin{figure}[ht!]
	\centering
	\begin{tikzpicture}[line cap=round,line join=round,>=triangle 45,x=1.3cm,y=1.3cm]
	\draw[line width=1pt] 
	(2.5,1.5)--(0,0)
	(2.5,1.5)--(1,0)
	(2.5,1.5)--(5,0)
	(2.5,1.5)--(6,0)
	(3.5,1.5)--(0,0)
	(3.5,1.5)--(1,0)
	(3.5,1.5)--(5,0)
	(3.5,1.5)--(6,0)
	(2,-1)--(0,0)
	(2,-1)--(1,0)
	(2,-1)--(5,0)
	(2,-1)--(6,0)
	(3,-1)--(0,0)
	(3,-1)--(1,0)
	(3,-1)--(5,0)
	(3,-1)--(6,0)
	(4,-1)--(0,0)
	(4,-1)--(1,0)
	(4,-1)--(5,0)
	(4,-1)--(6,0)
	(2,-1)--(3,-1)--(4,-1)
	(0,0)--(1,0)
	(5,0)--(6,0);
	\draw[line width=1pt] (2,-1) .. controls (2,-1.5) and (4,-1.5) .. (4,-1);
	\draw[line width=1pt] (0,0) .. controls (0,-0.25) and (5,-0.25) .. (5,0);
    \draw[line width=1pt] (0,0) .. controls (0,-0.35) and (5,-0.35) .. (6,0);
	\draw (3,0) node {$\ldots$};
	\draw[rounded corners](-0.3,0.4)--(6.3,0.4)--(6.3,-1.7)--(-0.3,-1.7)--cycle;
	\draw (5.5,-1.4) node {$K_{k+3}$};
	\draw (2.5,1.8) node{$x$};
	\draw (3.5,1.8) node{$y$};
        \draw (0,0.2) node{$v_1$};
        \draw (0.95,0.2) node{$v_2$};
        \draw (5.1,0.2) node{$v_{k-1}$};
        \draw (6.1,0.2) node{$v_{k}$};
        \draw (1.9,-1.2) node{$u$};
        \draw (3,-1.2) node{$v$};
        \draw (4.1,-1.2) node{$w$};
	\begin{scriptsize}
	\draw [fill=white] 
	(0,0) circle (2pt)
	(1,0) circle (2pt)
	(5,0) circle (2pt)
	(6,0) circle (2pt)
	(2.5,1.5) circle (2pt)
	(3.5,1.5) circle (2pt)
	(2,-1) circle (2pt)
	(3,-1) circle (2pt)
	(4,-1) circle (2pt);
	\end{scriptsize}
	\end{tikzpicture}
	\caption{Graph $Y_k$.}
	\label{fig:Graph $Y_k$}
\end{figure}

\begin{proposition}
\label{prop:dual-Y_k}
If $k\geq 2$, then $\gpd(Y_k) = 5$ and $\gpd(Y_k - x) = k+3$.
\end{proposition}

\proof
Let $\{v_1,\ldots,v_{k}, u,v,w\}$ be the set of vertices of the induced complete graph $K_{k+3}$ in $Y_k$ as indicated in Fig.~\ref{fig:Graph $Y_k$}. Then $V(Y_k) = \{v_1,\ldots,v_{k}, u,v,w,x,y\}$ and $\deg_{Y_k}(u) = \deg_{Y_k}(v) = \deg_{Y_k}(w) = k+2$. It is straightforward to check that $\{x,y, u,v,w\}$ is a general position set of $Y_k$. Hence $\{x,y, u,v,w\}$ is also a dual general position set of $Y_k$ by Theorem~\ref{thm:all}(iii). Therefore, $\gpd(Y_k)\geq 5$.

To prove that $\gpd(Y_k)\leq 5$, consider an arbitrary $\gpd$-set $X$ of $Y_k$.
Suppose to the contrary that $v_i\in X$ for some $i\in [k]$. Then exactly one of the vertices $x$ and $y$ must belong to $X$, for otherwise $v_i$ would lie on a shortest $x,y$-path, which is a contradiction.  Without loss of generality assume that $x\in X$ and $y\notin X$. We claim that there is no vertex from $\{u,v,w\}$ which is contained in $X$. If $z\in \{u,v,w\}\cap X$, then $v_i$ lies on a shortest $x,z$-path, but this contradicts the fact that $X$ is a $\gpd$-set of $Y_k$. If a vertex of $\{u,v,w\}$ does not contain in $X$, say $u$, then $v_i$ lies on a shortest $y,u$-path, which implies that $y$ and $u$ are not $X$-positionable. This contradiction yields $v_i\notin X$ for each $i\in [k]$. Therefore, $|X|\leq 5$, and we conclude that $\gpd(Y_k) = 5$.

If $x$ is removed from $Y_k$, by Theorem~\ref{thm:gpsets}, the set $\{v_1,\ldots,v_k,u,v,w\}$ is a largest general position set of $Y_k - x$. Using Theorem~\ref{thm:all}(iii), we get that 
$\{v_1,\ldots,v_k,u,v,w\}$ is also a dual general position set of $Y_k - x$. Note that $\gpd(Y_k-x)\leq \gp(Y_k-x) = k+3$. We can conclude that $\gpd(Y_k-x) = k+3$.
\qed

\section{Edge removal}
\label{sec:edge-removal}

Impact of removing an edge on the general position number of a graph was investigated in~\cite{dokyeesun-2025}. In this section, we complement this by examining the effect of edge removal on the other three general position invariants. 

We begin with the total general position number, for which the following notation is needed. If $G$ is a graph and $e=uv\in E(G)$, then let $S(G)_{e}$ denote the set of simplicial vertices in $G$ which are adjacent to both $u$ and $v$. 

\begin{proposition}
If $e$ is an edge of a graph $G$, then 
$$\gpt(G) -|S(G)_{e}| \le \gpt(G-e) \le \gpt(G) +2\,.$$
Moreover, the bounds are sharp. 
\end{proposition}

\proof
Let $e = uv$. If $x$ is a simplicial vertex which is adjacent to at most one of $u$ and $v$, then $x$ remains simplicial in $G-e$. Thus $\gpt(G-e) = s(G-e) \ge s(G) - |S(G)_e| = \gpt(G) - |S(G)_e|$. If $x\ne u,v$ is a simplicial vertex of $G-e$, it remains simplicial in $G$. Hence $\gpt(G) \ge \gpt(G-e) - 2$. 

For the sharpness, note first that $s(K_n) = n$ and $s(K_n-e) = 2$, hence the lower bound is sharp. For the upper bound, let $X_n$, $n\ge 2$, be the graph obtained from the disjoint union of two copies of $K_n$, say $ K$ and $K'$ and a vertex $x$, by adding an edge $e = uv$ between $K$ and $K'$, and joining $x$ to a vertex $u'$ in $K$ and a vertex $v'$ in $K'$, where $\{u,v\} \cap \{u',v'\} = \emptyset$. Then we infer that $\gpt(X_n)=2(n-2)$ and $\gpt(X_n-e)=2(n-1)$.
\qed

For the outer general position number, the following holds. 

\begin{theorem}
\label{thm:edge-removal-outer}
If $e$ is an edge in graph $G$, then
$$ \frac{\gpo(G)}{2}\le \gpo(G-e)\leq\ 2\gpo(G)\,.$$
Moreover, both bounds are sharp.
\end{theorem}

\proof
Let $e=uv$ and let $X$ be a $\gpo$-set of $G$. Then consider the sets of vertices 
\begin{align*}
 X_{uv} & = \{w\in X:\ d_G(u,w) < d_G(v,w) \}, \\
X_{vu} & = \{w\in X:\ d_G(v,w) < d_G(u,w) \}, \\
_vX_u & = \{w\in X:\ d_G(u,w) = d_G(v,w) \}\,.
\end{align*}
Clearly, $X = X_{uv} \cup X_{vu} \cup\, _vX_u$. Let 
$$X_u = X_{uv} \cup\, _vX_u\quad {\rm and}\quad X_v = X_{vu} \cup\, _vX_u\,.$$ 
Then we recall from the proof of~\cite[Theorem 6.2]{dokyeesun-2025} that $X_u$ and $X_v$ are general position sets of $G-e$. Moreover, by a parallel argument we prove that $X_u$ and $X_v$ remain to be MMD in $G-e$ provided they were MMD in $G$.  Since $|X|/2 \le \max\{ |X_u|, |X_v|\}$, Theorem~\ref{thm:all}(ii) yields the lower bound. The argument for the upper bounds proceeds along parallel lines. 

To demonstrate the sharpness of the lower bound, consider first graphs $Y'_n$, $n\ge 3$, constructed as follows. First, take two disjoint copies of $K_n$, say $K$ and $K'$ and add an edge between them, say $e=uv$, where $u\in K$ and $v\in K'$. Second, add two new vertices $u'$ and $v'$, the edge $u'v'$, all the edges between $u'$ and $V(K)\setminus \{u\}$, and all the edges between $v'$ and $V(K')\setminus \{v\}$, see Fig.~\ref{fig:Y'_4}.

\begin{figure}[ht!]
\begin{center}
\begin{tikzpicture}[scale=0.8,style=thick,x=1cm,y=1cm]
\def\vr{3pt}

\begin{scope}[xshift=0cm, yshift=0cm] 
\coordinate(x1) at (0,0);
\coordinate(x2) at (-1,2);
\coordinate(x3) at (0,2.5);
\coordinate(x4) at (1,2);
\coordinate(x5) at (0,4);
\coordinate(x6) at (3,0);
\coordinate(x7) at (2,2);
\coordinate(x8) at (3,2.5);
\coordinate(x9) at (4,2);
\coordinate(x10) at (3,4);
\draw (x1) -- (x2) -- (x3) -- (x4) -- (x5) -- (x3) -- (x1);
\draw (x6) -- (x7) -- (x8) -- (x9) -- (x10) -- (x8) -- (x6);
\draw (x5) -- (x10);
\draw (x2) -- (x4);
\draw (x7) -- (x9);
\draw (x1) -- (x6);
\draw (x2) -- (x5);
\draw (x7) -- (x10);
\draw (x1) -- (x4);
\draw (x6) -- (x9);
\foreach \i in {1,...,10} 
{
\draw(x\i)[fill=white] circle(\vr);
}
\foreach \i in {2,3,4,7,8,9} 
{
\draw(x\i)[fill=black] circle(\vr);
}
 circle(\vr);
\node at (-0.4,0) {$u'$}; 
\node at (-0.4,4) {$u$}; 
\node at (3.4,0) {$v'$}; 
\node at (3.4,4) {$v$}; 
\node at (1.5,4.3) {$e$}; 
\end{scope}

\begin{scope}[xshift=8cm, yshift=0cm] 
\coordinate(x1) at (0,0);
\coordinate(x2) at (-1,2);
\coordinate(x3) at (0,2.5);
\coordinate(x4) at (1,2);
\coordinate(x5) at (0,4);
\coordinate(x6) at (3,0);
\coordinate(x7) at (2,2);
\coordinate(x8) at (3,2.5);
\coordinate(x9) at (4,2);
\coordinate(x10) at (3,4);
\draw (x1) -- (x2) -- (x3) -- (x4) -- (x5) -- (x3) -- (x1);
\draw (x6) -- (x7) -- (x8) -- (x9) -- (x10) -- (x8) -- (x6);
\draw (x2) -- (x4);
\draw (x7) -- (x9);
\draw (x1) -- (x6);
\draw (x2) -- (x5);
\draw (x7) -- (x10);
\draw (x1) -- (x4);
\draw (x6) -- (x9);
\foreach \i in {1,...,10} 
{
\draw(x\i)[fill=white] circle(\vr);
}
\foreach \i in {2,3,4} 
{
\draw(x\i)[fill=black] circle(\vr);
}
 circle(\vr);
\end{scope}

\end{tikzpicture}
\caption{Graphs $Y'_4$ and $Y'_4-e$, and their $\gpo$-sets}
\label{fig:Y'_4}
\end{center}
\end{figure}

We can observe that $\gpo(Y'_n) = 2(n-1)$ and that $\gpo(Y'_n-e) = n-1$. This demonstrates the sharpness of the lower bound. 

For the sharpness of the upper bound, let $Z_n$, $n\ge 2$, be the graph obtained from two disjoint copies of $K_{2,n}$ by adding an edge $e=uv$ between them, where $u$ and $v$ are vertices of degree $n$ in different copies of $K_{2,n}$. Then $\gpo(Z_n) = n$ and $\gpo(Z_n-e) = 2n$. 
\qed

The sharpness of the upper bound in Theorem~\ref{thm:edge-removal-outer} is demonstrated by
graphs that contain bridges. Consider next the graphs $Z_n$, $n\ge 3$, constructed as follows. Take two disjoint complete graphs $K_n$, add an edge $e$ between them, and add a path of length $3$ between them disjoint from $e$. Then we can infer that $\gpo(Z_n) = 2n-4$ and that $\gpo(Z_n-e) = 2n-2$. In the case $n=3$, we thus also have a graph without bridges for which the upper bound in Theorem~\ref{thm:edge-removal-outer} is sharp.  It remains to be seen whether the upper bound in Theorem~\ref{thm:edge-removal-outer} is sharp also on graphs without bridges with an arbitrary large other general position number. 

To conclude the paper, we consider the dual general position number under edge removal for which we have the following. 

\begin{theorem}
\label{thm:gpd-gpd(G-e)} 
The difference $\gpd(G) - \gpd(G-e)$ can be arbitrarily large.
\end{theorem}

\proof
Consider the graphs $H_n$, $n\geq 1$, obtained from $n$ disjoint triangles and one $C_6$, by selecting an edge in each of these $n+1$ graphs and identifying them into a single edge. See Fig.~\ref{fig:H3} from which the vertex labelling of $H_n$ should be clear. 

\begin{figure}[ht!]
\begin{center}
\begin{tikzpicture}[scale=1.0,style=thick]
\tikzstyle{every node}=[draw=none,fill=none]
\def\vr{3pt} 

\path (0,0) coordinate (x);
\path (6,0) coordinate (y);
\path (3,1) coordinate (w1);
\path (3,2) coordinate (w2);
\path (3,4) coordinate (w3);
\path (1,-1) coordinate (v1);
\path (2.33,-1) coordinate (v2);
\path (3.66,-1) coordinate (v3);
\path (5.0,-1) coordinate (v4);
\draw (x)--(y);
\draw (x)--(w1)--(y);
\draw (x)--(w2)--(y);
\draw (x)--(w3)--(y);
\draw (x)--(v1)--(v2)--(v3)--(v4)--(y);
\draw (x)  [fill=white] circle (\vr);
\draw (y)  [fill=white] circle (\vr);
\draw (w1)  [fill=white] circle (\vr);
\draw (w2)  [fill=white] circle (\vr);
\draw (w3)  [fill=white] circle (\vr);
\draw (v1)  [fill=white] circle (\vr);
\draw (v2)  [fill=white] circle (\vr);
\draw (v3)  [fill=white] circle (\vr);
\draw (v4)  [fill=white] circle (\vr);
\draw (x)++(-0.3,0) node {$x$};
\draw (y)++(0.3,0) node {$y$};
\draw (w1)++(0.0,0.3) node {$w_1$};
\draw (w2)++(0.0,0.3) node {$w_2$};
\draw (w3)++(0.0,0.5) node {$w_n$};
\draw (w2)++(0.0,1.2) node {$\vdots$};
\draw (v1)++(0.0,-0.4) node {$v_1$};
\draw (v2)++(0.0,-0.4) node {$v_2$};
\draw (v2)++(0.66,0.8) node {$e$};
\draw (v3)++(0.0,-0.4) node {$v_3$};
\draw (v4)++(0.0,-0.4) node {$v_4$};
\end{tikzpicture}
\end{center}
\caption{The graph $H_n$}
\label{fig:H3}
\end{figure}

Let $X = \{w_1,\dots, w_n\}$. It is straightforward to check that $X$ is a general position set of $H_n$. Since there is a unique shortest path between $x$ and $y$ that does not contain any vertex from $X$, it follows that $G-X$ is convex. In addition, since $\{x,y,v_1,v_2,v_3,v_4\}$ induces a convex $C_6$, using Theorem~\ref{thm:all}(iii) we infer that none of these vertices lies in a dual general position set. We can conclude that $X$ is the largest dual general position set of $H_n$, therefore $\gpd(H_n) = n$. On the other hand, removing the edge $e = xy$, every edge of $H_n-e$ is the middle edge of some isometric path $P_4$. In view of~\cite[Proposition 3.3]{tian-2025} we have $\gpd(H_n-e) = 0$.
\qed

\section*{Acknowledgements}

This work was supported by the Slovenian Research and Innovation Agency (ARIS) under the grants P1-0297, N1-0355, N1-0285, N1-0431, and J1-70045.

\section*{Data Availability Statement}
 
Data sharing is not applicable to this article as no new data were created or analyzed in this study.

\section*{Conflict of Interest Statement}
 
Sandi Klav\v{z}ar is an Associate Editor of the Discrete Applied Mathematics journal and was not involved in the review and decision-making process of this article. In addition, the authors declare no other conflict of interest.

\baselineskip14pt



\begin{thebibliography}{99}

\bibitem{anand-2019}
B.S.~Anand, U.~Chandran S.V., M.~Changat, S.~Klav\v{z}ar, E.J.~Thomas,
Characterization of general position sets and its applications to cographs and bipartite graphs,
Appl.\ Math.\ Comput.\ 359 (2019) 84--89.

\bibitem{Araujo-2025}
J.~Araujo, M.C.~Dourado, F.~Protti, R.~Sampaio, 
The iteration time and the general position number in graph convexities, 
Appl.\ Math.\ Comput.\ 487 (2025) Paper 129084. 

\bibitem{survey}
U.~Chandran S.V., S.~Klav\v{z}ar, J.~Tuite, 
The general position problem: a survey,
\url{arXiv:2501.19385} (2025).

\bibitem{chandran-2016}
U.~Chandran S.V., G.J.~Parthasarathy,
The geodesic irredundant sets in graphs,
Int.\ J.\ Math.\ Combin.\ 4 (2016) 135--143.	

\bibitem{DiStefano-2025}
G.~Di Stefano, S.~Klav\v zar, A.~Krishnakumar, J.~Tuite, I.G.~Yero, 
Lower general position sets in graphs,
Discuss.\ Math.\ Graph Theory 45 (2025) 509--531. 

\bibitem{dokyeesun-2026}
P.~Dokyeesun, S.~Klav\v{z}ar, D.~Kuziak, J.~Tian,
General position problems in strong and lexicographic products of graphs,
Comput.\ Appl.\ Math.\ 45 (2026) Paper 97.

\bibitem{dokyeesun-2025}
P.~Dokyeesun, S.~Klav\v{z}ar, J.~Tian,
The general position number under vertex and edge removal, 
Quaest.\ Math.\ 48 (2025) 1277--1290. 

\bibitem{irsic-2024}
V.~Ir\v si\v c, S.~Klav\v zar, G.~Rus, J.~Tuite, 
General position polynomials,
Results Math.\ 79 (2024) Paper 110.

\bibitem{klavzar-2026}
S.~Klav{\v{z}}ar, A.~Krishnakumar, D.~Kuziak, E.~Shallcross, J.~Tuite, I.G.~Yero, 
Moving through {Cartesian} products, coronas and joins in general position,
Discrete Appl.\ Math.\ 379 (2026) 768--780. 

\bibitem{Korner-1995}
J.~K\"orner,
On the extremal combinatorics of the Hamming space,
J.\ Comb.\ Theory Ser.\ A 71 (1995) 112--126.

\bibitem{KorzeVesel-2025}
D.~Kor\v{z}e, A.~Vesel, 
Mutual-visibility and general position sets in {S}ierpi\'nski triangle graphs,
Bull.\ Malays.\ Math.\ Sci.\ Soc.\ 48 (2025) Paper 106.

\bibitem{Kruft-2025}
E.~Kruft Welton, S.~Khudairi, J.~Tuite,
Lower general position in Cartesian products,
Commun.\ Comb.\ Optim.\  10 (2025) 110--125.

\bibitem{manuel-2018a}
P.~Manuel, S.~Klav\v{z}ar,
A general position problem in graph theory,
Bull.\ Aust.\ Math.\ Sci.\ Soc.\ 98 (2018) 177--187.

\bibitem{oellermann-2007}
O.R.~Oellermann, J.~Peters-Fransen,
The strong metric dimension of graphs and digraphs,
Discrete Appl.\ Math.\ 155 (2007) 356--364.

\bibitem{Roy-2025}
D.~Roy, S.~Klav\v{z}ar, A.S.~Lakshmanan,
Mutual-visibility and general position in double graphs and in Mycielskians,
Appl.\ Math.\ Comput.\ 488 (2025) Paper 129131.

\bibitem{ThaChaTuiThoSteErs-2024} 
M.~Thankachy, U.~Chandran S.V., J.~Tuite, E.~Thomas, G.~Di Stefano, G.~Erskine, 
On the vertex position number of graphs, 
Discuss.\ Math.\ Graph Theory 44 (2024) 1169--1188.  

\bibitem{Thomas-2024a}
E.J.~Thomas, U.~Chandran S.V., J.~Tuite, G.~Di Stefano,
On the general position number of Mycielskian graphs,
Discrete Appl.\ Math.\ 353 (2024) 29--43. 

\bibitem{tian-2025}
J.~Tian, S. Klav\v{z}ar,
Variety of general position problems in graphs,
Bull.\ Malays.\ Math.\ Sci.\ Soc. 48 (2025) Paper 5. 

\bibitem{tian-2023}
J.~Tian, K.~Xu, D.~Chao,
On the general position numbers of maximal outerplane graphs,
Bull.\ Malays.\ Math.\ Sci.\ Soc.\ 46 (2023) Paper 198.

\bibitem{welton-2025}
E.K.~Welton, S.~Khudairi, J.~Tuite, 
Lower general position in Cartesian products,
Commun.\ Comb.\ Optim.\ 10 (2025) 110--125.

\end{thebibliography}
\end{document}